\documentclass[a4paper,10pt]{article}
\usepackage[utf8x]{inputenc}

\usepackage{amssymb}

\title{Special Bohr - Sommerfeld geometry.}
\author{{\bf Nikolay A. Tyurin}\\
BLTPh (Dubna)\footnote{Joint Institute for Nuclear research, Joliot - Curie, 6, Moscow region, Dubna 141980, Russia}\\ {\it and } \\
NRU HSE (Moscow)\footnote{AG Laboratory,
 HSE, 7 Vavilova str., Moscow, Russia, 117312}}

\begin{document}

\maketitle

\begin{abstract} We introduce a new notion --- special lagrangin submanifolds,
which satisfy the Bohr - Sommerfeld condition --- for algebraic varieties. We show that this leads to the construction
of finite dimensional moduli space of special Bohr - Sommerfeld lagrangian submanifolds with respect to any ample linear bundle.
The construction can be used in the studies of Mirror Symmetry.

\end{abstract}

The essence of Mirror Symmetry in the broadest context was expressed by Yu. I. Manin as ``duality between 
symplectic geometry and complex Kahler geometry'' (see [1]). 
Two algebraic Kahler manifolds $M, W$ are understood as ``mirror partners'', if certain derived objects,
constructed in the frameworks of algebraic and symplectic geometries of $M, W$ are cross equivalent: for example
in Homological Mirror Symmetry due to M. Kontsevich (see [2]) the derived category of coherent sheaves 
on $M$ must be equivalent to the Fukaya - Floer category of $W$ and vice versa.  

A.N. Tyurin, who spent many years studing stable vector bundles, suggested more geometrical correspondence:
certain duality between vector budnles and lagrangian submanifolds (see, f.e. [3]). Namely
for a pair of threefolds  $M, W$ even cohomologies represent the Chern classes of vector bundles, combined into 
finite dimensionalmoduli spaces of stable vector bundles, and the middle odd cohomology can be realized 
by lagrangian submanifolds, should be combined into finite dimensional moduli spaces; and then on the comparing
of these moduli spaces one could define the duality, which should present the essence of Mirror Symmetry.
The main problem arises in this way  is in   ``infinitness'' of lagrangian geometry in contrast to
algebraic geometry.

The problem can be solved if one introduces certain speciality condition on lagrangian submanifolds: 
realizing the ideology of calibrated lagrangian cycles J. MacLean and N. Hitchin (see [4])
proposed a special condition on lagrangian submanifolds in Calabi - Yau varieties which led to
finite dimensional moduli spaces. Briefly, any Calabi - Yau variety by the very definition is endowed by top holomorphic 
non vanishing form $\theta \in \Omega^{3,0}$, and its restriction to any lagrangian (w.r.t. the Kahler form) $S$
is non vanishing as well; therefore the condition $\theta|_S = \psi d\mu(g)$ 
gives a correclty defined complex function $\psi_S: S \to \mathbb{C}^*$. One says that  $S$ is special (or SpLAG) iff $\psi$
has constant argument. Local deformations of the SpLAG submanifolds are finite dimensional and unobstructed,
so the moduli space of special lagrangian submanifolds in a Calabi - Yau threefold has dimension $b_1(S)$ (details
can be found in [4]).

 The introduction of SpLAG geometry led to the realization of Mirror Symmetry as T - duality:
 the famous  SYZ - conjecture for Calabi - Yau threefolds (see [5]) explains Mirror Symmetry in terms
 of fibrations on SpLAG tori. According to this conjecture, kahlerian Calabi - Yau threefolds
 can be fibered on special lagrangian tori, parameterized by certain three dimenisonal base $B$ (note that $b_1(T^3) = 3$),
 and the mirror partner is given by the fiberation on dual tori over the same base $B$ (see [5]).

Unfortunately the existence problem for such a special lagrangian fibrations on Calabi - Yau threefolds is still open
despite of crucial attempts to proof, and after years of high level popularity of SpLAG geometry the community 
of "mirror symmetrists" has changed
the focus to Kontsevich's homological approach. D. Auroux in [6] revisited the subject, extending the notion of SpLAG submanifolds
to the case of open Calabi -Yau varieties, given by cut of a divisor from the anticanonical linear system on Fano varieties.
In this approach the modified speciality condition turns to be relative --- for a given Fano variety (or, more rigourisly,
for a variety, whose anticanonical bundle admits holomorphic sections) this condition depends on the choice
of the divisor and a lagrangian submanifold must lie in the complement of this divisor to be special, 
and since this complement can be understood as an open Calabi - Yau manifold, the Auroux's approach can be seen as 
a modification of SpLAG geometry. The simplest example of special lagrangian fibration of an open Calabi -Yau is given
by toric geometry: removing three lines $l_i = \{ z_i = 0 \}$ from the projective plane $\mathbb{C} \mathbb{P}^2$
one gets a special lagrangian fibration by Clifford tori. In [6] it was constructed a special lagrangian fibration
for the complement of a reducible cubic equals to the union of nondegenerated quadric and projective line, and then
it was conjectured that such a fibration exists for the complement of a smooth cubic curve in $\mathbb{C} \mathbb{P}^2$,
but this conjecture is still open. Except a single example from [7] there is no activity in this way.

In [8] one developed ideas of lagrangian approach to geometric quantization, namely one studied the moduli space
of lagrangian submanifolds which satisfy the Bohr - Sommerfeld condition. In a previous paper A.N. Tyurin  observed 
that the Bohr - Sommerfeld condition is in a sense "transversal" to the SpLAG conidtion in the Calabi - Yau case,
which leads to possible definitions of certain finite invariants for Calabi - Yau threefolds, mirror to the Cassons
invariants (see [3]). We develop this idea and introduce a new speciality condition for Bohr - Sommerfeld lagrangian
submanifolds only, and this new notion of special Bohr - Sommerfeld submanifolds opens a chain of intersting
observations, collected by Special Bohr - Sommerfeld geometry (SBS geometry for short).

Let $(M, \omega)$ be a compact simply connected symplectic manifold, satisfied the <bohr - Sommerfeld condition of symplectic
manifolds: the cohomology class  $[\omega] \in H^2(M, \mathbb{R})$ is integer. Then fix the so called prequantization
data: a linear bundle $L \to M$ endowed by a hermitian structure and a hermitian connection $a$ on it such that
the curvature form  $F_a = 2 \pi i \omega$. This condition uniquelly defines the connection up to gauge
transformations. A lagrangian submanifold $S \subset M$ satisfies the Bohr - Sommerfeld condition (BS for short) iff
the restriction $(L, a)|_{S}$ admits a covariantlyconstant section $\sigma_S$. 

Let  $s \in \Gamma(M, L)$ be any smooth section of  $L$.

{\bf Definition.} {\it We call a BS- lagrangian submanifold $S$ special with respect to section $s$, 
iff  $s|_S$ nowhere vanishes and the proportionality cofficient  $\alpha(s, S)$, defined by the equality 
$s|_S = \alpha(S,s) \sigma_S$, has constant argument.}

Since the definition doesn't depend on the choice of $\sigma_S$ and of rescaling of $s$, it induces an "incidence 
cycle"in the direct product
$$
{\cal U}_{SBS} \subset \mathbb{P}(\Gamma(M, L)) \times {\cal B}_S,
$$
where the last symbol denotes the moduli space of Bohr - Sommerfeld lagrangian cycles of fixed topological type (see [8]).
Namely, a pair  $([s], S)$ belongs to ${\cal U}_{SBS}$ iff  $S$ is SBS w.r.t section $s$, if $s$ represents the class
$[s]$ in the projectivized space. Naturally there are two projections $p_1, p_2$ to the first and to the second direct
summands.

``Finitness'' of the set of SBS submanifolds is reflected by the following fact:

{\bf Theorem A.} {\it The projection $p_1: {\cal U}_{SBS} \to \mathbb{P}(\Gamma (M, L))$ has discrete fibers 
over the image $\rm{Im} p_1 \subset \mathbb{P}(\Gamma (M, L))$.}

On the other hand, one has

{\bf Theorem B.}{\it The image $\rm{Im} p_1 \subset \mathbb{P} (\Gamma (M, L))$ is an open subset in the projective space.}

Together they give

{\bf Corollary A.} {\it The space ${\cal U}_{SBS}$ admits a Kahler structure.}

Note that in general this fact leads to possible applications in Geometric Quantization.
However we are strongly interested in the specific case: suppose that our symplectic manifold $(M, \omega)$ 
admits an integrable complex structure $I$, compatible with $\omega$. This means that $M$ is algebraic
variety with principal polarization, defined by a holomorphic line bundle $L$.
Then we get a finite dimensional subspace $\mathbb{P}(H^0(M_I, L)) \subset \mathbb{P}(\Gamma(M, L))$ 
formed by classes of holomorphic sections and therefore a reduced "incidence cycle"
$$
{\cal M}_{SBS} \subset \mathbb{P}(H^0(M_I, L)) \times {\cal B}_S
$$
together with two natural projections to the direct summands which we again denote as $p_1, p_2$. Then we have

{\bf Corollary B.} {\it The moduli space ${\cal M}_{SBS}$ is finite dimensional possible singular Kahler variety.}

Thus for any compact simply connected algebraic variety $X$ one can construct the following family of
 moduli spaces: for each very ample $L \to X$ the corresponding complete linear system defines the embedding of 
 $X$ to the projective space, dual to $\mathbb{P} H^0(X, L)$; we can lift to $X$ the standard Kahler form
 from the projective space and consider it as a symplectic form on $X$ --- and this is exactly the case
 where SBS geometry can be switched on. For fixed topological type of $S$ and the class  $[S] \in H_n(X, \mathbb{Z})$
in the middle cohomology of  $X$ we get the corresponding moduli space ${\cal M}_{SBS}$ marked by the data
$$
{\cal M}_{SBS} = {\cal M}_{SBS} (S, [S], c_1(L)), \quad c_1 (L) \in H^2(X, \mathbb{Z}).
$$  
Note however that even for close bundles these moduli spaces can be drastically different, and for the same bundle 
$L$ and the same class $[S]$ but for different topological types of the moduli space can be different as well:
below we give the examples.

{\bf Example 1.} Consider as $M$ the simplest compact simply connected symplectic manifold ---
complex projective line $\mathbb{C} \mathbb{P}^1$ endowed by the standard Kahler structure.
If we take as $L$ the bundle  ${\cal O}(1)$, then the moduli space  ${\cal M}_{SBS} (S^1, 0, h)$ 
is empty since no smooth loops on $\mathbb{C} \mathbb{P}^1$ satisfy the Bohr - Sommerfeld condition (see [8]). 
But if we take as $L$ the bundle ${\cal O}(2)$, then the moduli space ${\cal M}_{SBS} (S^1, 0, 2h)$ 
is naturally isomorphic to  $\mathbb{C} \mathbb{P}^2 \backslash
Q$, where conic $Q$ is the image of $\mathbb{C} \mathbb{P}^1$ under the Veronese embedding (see Section 3).

{\bf Example 2.} For $M = \mathbb{C} \mathbb{P}^2$ with the standard Kahler structure if we take
${\cal O}(2)$ as $L$ then the moduli space ${cal M}_{SBS} (T^2, 0, 2h)$ is empty while the moduli space
${\cal M}_{SBS}(\mathbb{R}\mathbb{P}^2, 0, 2h)$ is nonempty. The last  fact can be seen from the following arguments:
fix coordinates $[z_0: z_1: z_2]$ compatible with the fixed Kahler structure and consider $S = \mathbb{R} \mathbb{P}^2 =
\{ z_i \in \mathbb{R} \}$. Then it is not hard to see that $S$ is SBS w.r.t. the holomorphic section with zeros
on the Fermat conic $Q = \{z_0^2 + z_1^2 + z_2^2 = 0\}$.

These results are derived in view of the following observation: the special Bohr - Sommerfeld condition
with respect to a holomorphic section is naturally related to the Morse theory of plurisubharmonic
functions on the complements to divisors. As we show below, if $s$ is a holomorphic section
with zeroset $D_s \subset M$ then one takes $\phi_s = - ln \vert s \vert^2$ which is plurisubharmonic on
$M \backslash D_s$ and then SBS condition for lagrangian submanifold $S \subset M$
reads as $grad \phi_s \subset TS$ at each point of $S$. Equivalently $S$ is preserved by the gradient flow of 
$\phi_s$. In particular this means that if $\phi_s$ is Morse outside of $D_s$ then the number and
the type of critical points of $\phi_s$ dictate the possible types of lagrangian submanifolds.
These critical points must be critical points of a Morse function on $S$, therefore to have
a SBS torus of dimension 2  one must have at least 4 critical points of $\phi_s$ otuside of
$D_s$, which is non realistic for holomorphic sections of ${\cal O}(2)$; in contrast the set of critical
points for $\phi_s$ where $s$ corresponds to the Fermat conic is very big --- they form exactly
$\mathbb{R}\mathbb{P}^2$, and it is SBS due to the reformulated in Kahler terms SBS condition.

These arguments lead to

{\bf Theorem C.} {\it For a generic holomorphic section $s \in H^0(M_I, L)$ the set of special Bohr - Sommerfeld
lagrangian submanifolds is finite.}

We expect that it is true for any holomoprhic section, and that one has in general certain ramified covering
$$
{\cal M}|_{SBS} \to \mathbb{P} (H^0(M_I, L)),
$$
so should get a reach geometrical picture combined lagrangian and Kahler geometries of $M$.

The present paper is the first step and draft, so we focus on ideas and leave aside some technical details
and computations in proofs. 

{\bf Acknowledgments.} The defintion of special Bohr - Sommerfeld submanifold
arose  in the discussion with Andery Losev, and the author would like to thank
him first of all. I'm gratefull to Andrey Shafarevich for the dicussion of general
constructions and Vsevold Schevchishin who pointed out the relations to
the theory of Kahler potential.  I would like to thank Anastassia Tyurina for the computation
in the toy example $M = \mathbb{C} \mathbb{P}^1$.

The article was prepared within the framework of a subsidy granted to the HSE by the
Government of the Russian Federation for the implementation of the Global Competitiveness Program.

\section{Special Bohr - Sommerfeld submanifolds}

Consider a compact simply connected symlectic manifold $(M, \omega)$ of real dimension $2n$ 
such that the symplectic form $\omega$ has integer class in the de Rham cohomology of $M$: $[\omega] \in H^2(M, \mathbb{Z})$.
In this case one says that the symplecticmanifold $(M, \omega)$ satisfies the Bohr - Sommerfeld condition for
manifolds. Fix the prequantization data $(L, a)$: linear hermitian bundle $L$ whose first Chern class is
$c_1 (L) = [\omega] \in H^2(M, \mathbb{Z})$ and a hermitian connection  $a \in {\cal A}_h(L)$,
with the curvature form is $F_a = 2 \pi i \omega$. In the simply connected case this condition defines $a$ 
uniquelly up to gauge transformations. Recall that the quadruple $(M, \omega, L, a)$ is the input data
for the ALAG programme, see [8]. 

A submanifold $S \subset M$ of real dimension $n$ is lagrangian iff the restriction $\omega|_S \equiv 0$ is trivial.
Therefore for any lagrangian submanifold the restriction $(L, a)|_S$ is topologically trivial
bundle with a flat connection. In this paper we mostly consider the case of compact orientable
lagrangian submanifolds although certain results can be extended to more general cases.

{\bf Definition 1.} {\it  Lagrangian submanifold $S$ satisfies the Bohr - Sommerfeld condition iff
the restriction  $(L, a)|_S$ admits a covariantly constant section $\sigma_S \in \Gamma(S, L|_S)$. }

In what follows we will abbreviate it as BS - condition. 

It is not hard to see that  BS - condition doesn't depend on the choice of $a$ in the equivalence class
up to gauge transformations. Indeed, {\bf Definition 1} is equivalent to the following condition:
for any loop  $\gamma \subset S$ and any disc $D \subset M$ such that $\partial D = \gamma$
the symplectic area of disc $D$ is integer: $\int_D \omega \in \mathbb{Z}$. 
At the same time for a fixed $a$ the covariantly constant section $\sigma_S$ is defined up to  $\mathbb{C}^*$.

Fix a smooth section $s \in \Gamma (M, L)$ of the prequantization bundle $L \to M$. Then for a BS - submanifold
 $S \subset M$ we give the following 

{\bf Definition 2.} {\it We say that BS- submanifold $S$ is special w.r.t section $s$, iff $s|_S$ doesn't vanish
on $S$ and the proportionality coefficient $\alpha(s, S)$, defined by the equality $s|_S = \alpha(S,s) \sigma_S$, 
has constant argument. In other words $s|_S =
f e^{ic} \sigma_S$, where $c$ is a real constant and $f$ --- a real positive function.}

 For short we will call such an $S$ as $s$ - SBS - submanifold or just SBS - submanifold if
 $s$ is coming from the context. 

{\bf Remark.} Our SBS - condition presented above  is essentially different from
the speciality conditions given by N. Hitchin in [4] and by D. Auroux in [6], despite of the fact
that the second one depends on the section as in our case. Our SBS - condition is based on 
the Bohr - Sommerfeld condition and forother lagrangian submanifold it can't be
either generalized nor reformulated. At the same time certain weak relation takes place. Namely 
 consider the case when $K_M = k [\omega]$, where $K_M = det ((T^*M)^{1,0})$ is the canonical class 
of a Kahler manifold $M$,  and connection $a$ is defined by the condition $F_a = 2 \pi i k \omega$. 
Then both in the Calaby - Yau case in  [4] ($k =0$) and in the Fano case in  [6] ($k<0$) one can take 
the determinant Levi - Civita connection as $a$ (since the Kahler metric is the Kahler - Einstein), 
and if the mean curvature of our lagrangian submanifold $S$ identically vanishes then the covarinatly constant section
$\sigma_S$ is up to multiple the volume form of the restricted Kahler metric to $S$. Therefore
in the both cases {\it minimal} and special in the sense of Hitchin or Auroux lagrangian submanifolds
are SBS - submanifolds with respect to the corresponding sections of the anticanonical bundle,
given by  top holomorphic form $\theta$ either on whole  $M$ as in [4], or on the complement 
to the corresponding divisor form the anticanonical system as in [6]. 
Note however that lagrangian submanifold $S$ is minimal only if it is Bohr - Sommerfeld
with respect to the determinant Levi - Civita connection, but BS - condition is not sufficient
in the case  so there are  SpLag - submanifolds and BS - lagrangian submanifolds which are not SBS. 
However we think that it is reasonable to use term "special" in our  Definition  1 above.

Recall that ALAG - programme from [8] supplies one in the situation described above with certain moduli space ${\cal B}_S$
of Bohr - Sommerfeld lagrangian submanifolds of fixed topological type (see [8]). Every such a moduli space
is defined and numerated by the following discrete data: by the corresponding class $[S] \in H_n(M, \mathbb{Z})$, 
realized by the lagrangian submanifolds, and by the topological type of $S$ (f.e. in the case  $n=2$
the lasr type is fixed by the genus $g(S)$); we indicate this dependence as ${\cal B}_S(S, [S])$. This moduli space
is an infinite dimensional smooth by Frechet real manifolds. At each point its tangent  space is modelled
by the real space $C^{\infty}(S, \mathbb{R})/const$ (the details can be found in  [8]).

It's clear that SBS - condition is stable with respect to rescaling of
smooth section $s$, hence this condition induces certain "incidence cycle" in the direct product
$\mathbb{P} \Gamma(M, L) \times {\cal B}_S$. Namely define a subset $\mathbb{P} \Gamma (M, L) \times {\cal B}_S \supset 
{\cal M}_S = \{ (p, S)\}$ by the condition that  BS - submanifold $S \subset M$
is $s$- SBS - submanifolds w.r.t the section $s \in \Gamma(M, L)$, representing the point $p \in \mathbb{P}\Gamma(M, L)$.

As usual, for our "incidence cycle" ${\cal U}_{SBS}$ one has two canonical projections
$$
p_1 : {\cal U}_{SBS} \to \mathbb{P} (\Gamma (M, L)), \quad p_2: {\cal U}_{SBS} \to {\cal B}_S
$$
to the direct summands of the ambient direct product. The main part of the present work
is to study the properties of these projections since the geometrical properties of the summands
 $\mathbb{P} (\Gamma (M, L))$ are ${\cal B}_S$ already known. Moreover both of them play essential roles
 in the Geometric Quantization constructions.
 
{\bf Digression: Geometric Quantization.} Infinite dimensional projective space $\mathbb{P}(\Gamma(M, L)$ ---
one of the principal objects in Geometric Quantization. Usually in this projective space one cuts subspaces
which represent quantum phase spaces of the qunatized system. However the second summand --- the moduli space
${\cal B}_S$ --- has been exploited in other approaches to the quantization problem for classical mechanical systems
named as Lagrangian Geometric Quantization. In this set up one exploited certain building over  ${\cal B}_S$,
given in  [8], when after "halfweighting" one constructs  ${\cal B}_S^{hw, r}$ --- the moduli space of halfweighted
Bohr - Sommerfeld lagrangian cycles of fixed topological type and volume. The aim of this building 
was a "complexification" of real moduli space ${\cal B}_S$. The moduli space ${\cal B}_S^{hw, r}$ 
was used in the construction named as ALG(a) - quantization, see [9]. But the original  aim wasn't reached in this way
since this building admits an almost Kahler structure with the constant Kahler angle which is not integrable.
But in any case our "incidence cycle" could play an important and interesting role: being universal
in the direct product it could helps to carry geometrical data from  $\mathbb{P} (\Gamma (M, L))$ to ${\cal B}_S$ 
and back. Thus it could give a connection between different approaches to Geometric Quantization
problem. Moreover, since (as we will see below) the projection $p_2$ is epimorphic and since ${\cal U}_{SBS}$ 
admits a Kahler structure lifted from the first summand, then our "incidence cycle"
can be regarded as a complexification of ${\cal B}_S$.

Defitinion 2 implies the following simple topological observation. In the direct product $\mathbb{P}(\Gamma (M, L)) \times {\cal B}_S$
take the determinantal subset $\Delta \subset  \mathbb{P}(\Gamma (M, L)) \times {\cal B}_S$ by the condition: 
section $s$ after the restrictionto  $S$ has zeros. Then the complement $\mathbb{P}(\Gamma (M, L)) \times {\cal B}_S
\backslash \Delta$ can be divided into connected components $K_i$, possible of infinite number. For each connected component 
$K_i$ one can define the cohomology class $m(K_i) \subset H^1(S, \mathbb{Z})$ by the  condition: for pair $(p, S) \subset K_i$
the resctriction of  $s$ to  $S$ gives the function $\alpha (s, S) \in C^{\infty}(S, \mathbb{C}^*$, so
$$
m(K_i) = \alpha(s, S)^* \mu
$$, where  $\mu$ is the generator in  $H^1(\mathbb{C}^*, \mathbb{Z})$. 
Inside of the connected component  the class $m(K_i)$ can't change being integer valued therefore the definitionof $m(K_i)$
doesn't depend on the choice of $(p, S) \in K_i$.

Then one has 

{\bf Proposition 1.} {\it The moduli space ${\cal U}_{SBS}$ has non trivial intersection with the component $K_i$ iff
$m(K_i) = 0$.}

The proof is  obvious. 

The next proposition is more interesting: it reflects the most important for us geometrical property of the first projection
 $p_1$.

{\bf Proposition 2.}  {\it For any smooth section  $s \in \Gamma (X, L)$ the set of $s$ - SBS submanifolds of
fixed topological type is discrete.}

{\bf Proof.} Suppose in contrary that there is one dimensional family $S_t, t \in [0, 1),$ 
of SBS - lagrangian submanifolds for a fixed smooth section $s_0$. Consider a Darboux - Weinstein
neighborhood ${\cal O}_0$ of the origin lagrangian submanifold $S_0$; then there is a sufficiently small segment
$[0; \epsilon]$ such that every $S_t \subset {\cal O}_0$ if $t \in [0; \epsilon]$ and consequently each $S_t$
is presented by the corresponding exact 1 - form $d \psi_t \in \Gamma (S_0, T^* S_0)$.
Hence each $S_t$  intersects $S_0$ at least at two points $p^+_t, p^-_t$ if $t \in [0; \epsilon]$ since every
function on a compact set must have at least two critical points.
Join these two points by pathes $\gamma_0, \gamma_t$ which lie on the submanifolds $S_0$  and $S_t$ respectively. 
Choice covariantlyconstant sections $\sigma_0 \in \Gamma (S_0, L|_{S_0})$ and $\sigma_t \in \Gamma (S_t, L|_{S_t}$
such that $\sigma_0(p^-t) = \sigma_t(p^-_t)$. Then it is clear that the phase difference at the second common point $p^+_t$ 
can be expressed via the symplectic area of the disc  $D_t$, bounded by  $\gamma_0$ and $\gamma_t$:
$$
\frac{\sigma_t(p_t^+)}{\sigma_0(p_t^+)} = {\rm exp} [2 \pi i \int_{D_t} \omega] = 2 \pi i (\psi_t(p_t^+) - \psi_t (p_t^-)).
$$
But both the BS - submanifolds $S_0$ and $S_t$ are special w.r.t the same global section  $s_0$, therefore the phase of $ \sigma_0$
and  $\sigma_t$ must coincide at both the points $p^{\pm}_t$, which is impossible for sufficiently small $t$ 
if $\psi_t(p_t^+) - \psi_t (p_t^-)$ is not zero. But if it is then the maximum and the minimum of our function
coincide so the function is constant and hence $S_t = S_0$.

{\bf Remark.} If one removes from the Definition 2 the non vanishing condition for the restriction
of $s$ to $S$ (we would call this "stability condition") then it is easy to construct 
a contineous deformation of SBS - submanifolds for the same section $s_0$. Indeed, as the simplest example 
we take  $\mathbb{C}$ with the standard symplectic form $\frac{1}{2 \pi i} dz \wedge d \bar z$,
then the prequantization bundle it trivial $L = C^{\infty}(\mathbb{C}, \mathbb{C})$. The prequantization
connection with respect to the natural trivialization is given by 1-form $\frac{1}{2 \pi} (z d \bar z - \bar z dz)$, 
and any real line of the form $z = ct$ is  BS - submanifold such that the correspondings covariantly
constant sections are given by constant functions. Holomorphic section $f(z) = z$ vanishes at $z=0$, 
and the pencil of real line $z = ct$ is a continous family of SBS - submanifolds, intersecting exactly at
the origin. The situation can be modified to the compact case if we complete  $\mathbb{C}$ to the projective 
line $\mathbb{C} \mathbb{P}^1$ and consider the family of meridians passing thorugh the Poles ---
all of them are SBS in this weaker sense (non stable) w.r.t to the section given by the homogenious polynomial
 $z_0 z_1$, for the prequantization bundle $L = {\cal O}(2)$.

Also note that $p_1$  is never epimorphic: one can take a smooth section $s \in \Gamma (M, L)$ with very big set of zeros
$(s)_0 = \{ x \in M | s(x) = 0 \} \subset M$ such that the complement $M \backslash (s)_0$ would not contain a disc $D$
with integer symplectic area (the argument doesn't work when $M$ admits family of shrinking lagrangian spheres).

Local computation ensures that

{\bf Proposition 3.} {\it Over a generic point $p \in Im p_1$ the differential of $d p_1$ is an isomorphism.}

Therefore we can state that

{\bf Theorem 1.} {\it The projection  $p_1: {\cal U}_{SBS} \to \mathbb{P}(\Gamma (M, L))$ has the structure of 
 covering  over the image $Im p_1 \subset \mathbb{P}(\Gamma (M, L))$.}
 
 Below we continue the studies of the first projection $p_1$.

\section{Speciality and calibration}

In this section we present an alternative description of
special Bohr - Sommerfeld in calibration terms. Recall that calibration in the sense of Harvey and Lawson [11] 
is given by the vanishing condition for a set of distingished form after the restriction to submanifolds. 
Below we show that SBS - condition is equivalent to vanishing of certain 1 - form, constructed in terms of
section $s$.

Consider $s \in \Gamma (M, L)$ and denote its zeroset as $D_s \subset M$. Then on the complement
$M \backslash D_s$ it is correctly defined the following complex 1 - form $\rho_s \in \Omega^1_M \otimes \mathbb{C}$
$$
\rho_s = \frac{\nabla_a s}{s} = \frac{<\nabla_a s, s>}{<s, s>},
$$
where $\nabla_a$ is the covariant derivative of the prequantization connection $a \in {\cal A}_h(L)$.

Indeed, $\nabla_a s \in \Gamma (L) \otimes \Omega^1_M$, and since $s$ on $M \backslash D_s$ doens't vanish we can
express  $\nabla_a s$ as a section of $L \otimes T^*M$ in the form $s \otimes \rho_s$ outside of
zeroset of $s$.

{\bf Proposition 4.}{\it Complex  1 - form  $\rho_s$ doesn't change under rescaling of $s$ by a constant.
The real part of $rho_s$ is exact. The differential of the imaginary part of $\rho_s$ equals $2 \pi \omega$.} 

First, rescale $s$ by a constant $c \in \mathbb{C}$ doesn't change the zeroset $D_s$, at the same time
$$
\frac{\nabla_a c s}{cs} = \frac{c \nabla_a s}{ c s} = \frac{\nabla_a s}{s}
$$
on $M \backslash D_s = M \backslash D_{cs}$. Therefore pair ($D_s \subset M$; 1 - form $\rho_s$) corresponds
to a point from $\mathbb{P}(\Gamma (M, L))$.

Second, the real part of $\rho_s$ can be found from the equality:
$$
d< s, s> = <\nabla_a s, s> + <s, \nabla_a s> = 2 \rm{Re} \rho_s <s, s>,
$$
since $a$ is hermitian. Therefore
$$
\rm{Re} \rho_s = 1/2 d (\rm{ln} \vert s \vert^2),
$$
so the part is exact on  $M \backslash D_s$.

Third, the calculations  of the imaginary part of $\rho_s$ can be doneas follows:
take on $M \backslash D_s$ (non hermitian) connection $a_s$ such that $\nabla_{a_s} s = 0$. 
Then the connection $a_s$ is trivial and defined by the trivialization given by
the section $s$ of the restricted bundle $L|_{M \backslash D_s}$. But the difference $\nabla_a - \nabla_{a_s}$ 
is exactly $\rho_s$, and the differential of $\rho_s$ must be equal to the difference of the curvature forms 
$F_a - F_{a_s} = d \rho_s$. The first one by the very definition equals to  $2 \pi i \omega$, 
and the second one is trivial. It follows $d \rho_s = 2 \pi i \omega$, so
$$
d (\rm{Re} \rho_s) = 0, \quad d (\rm{Im} \rho_s) = 2 \pi \omega,
$$
which ends the proof.  

{\bf Remark.} Proposition 4 implies that any section $s \in \Gamma (M, L)$ defines the following structure on
the open part $M \backslash D_s$ where $D_s$ is the zero set of $s$: a smooth function $\phi_s = - \rm{ln}\vert s \vert$
bounded below and a vector field $X_s$ given by $\frac{1}{2 \pi}\omega^{-1} (\rm{Im} \rho_s)$, which is a Liouville
vector field since the Lie derivative $Lie_{X_f} \omega = \omega$. Both $\phi_s$ and $X_s$ are ecnoded by our 1 -form
$\rho_s$. The structure given by $(\omega, X_s)$ is very well known as the Liouville structure, and
moreover if the vector field $X_s$ is {\it gradient - like} for the function $\phi_s$ then we get
a Weinstein structure on $M \backslash D_s$. It seems that it happens iff there exists a compatible
almost complex structure $J$ such that the section $s$ is pseudo holomorphic w.r.t. $J$. We will
exploit the way ``back from Weinstein to Stein'' built up  by Y. Elieasberg and K. Cielieback\footnote{see
the lecture course at www.mathematik.uni-muenchen.de/~kai/research/stein.pdf}
in the next section devoted to the algebraic case.

The point is that our  1- form $\rm{Im} \rho_s$ defines the requared calibration:

{\bf Theorem 2.} {\it Smooth orientable submanifold $S \subset M$ of dimension $n$ is  $s$ - special Bohr - Sommerfeld
lagrangian iff the restriction of 1 - form $\rm{Im} \rho_s$ to $S$ identically vanishes.}

Note that the condition $(\rm{Im} \rho_s)|_S \equiv 0$ implies that  $S \cap D_s = \emptyset$ 
since the form $\rho_s$ has pole along  $D_s$, therefore we don't mention the last condition
in the Theorem formulation.  

{\bf Proof.} Let $S$ be special w.r.t. $s$ lagrangian BS - submanifold. Then the restriction of $\rho_s$ on
$S$ can be calculated as follows: 
$$
\rho_s|_S = \frac{\nabla_a|_S (s|_S)}{s|_S} = \frac{\nabla_a f e^{ic} \sigma_S}{f e^{ic} \sigma_S} = \frac{df}{f} = d (\rm{ln} f),
$$
where $f$ is real positive function (see Definition 2). Thus $(\rm{Im} \rho_S)|_S \equiv 0$. 

Note that during the calculation above we establish the meaning of positive function $f$ ---
it is exactly $ \vert s \vert$, restricted to $S$.  

Now let certain smooth orientable half dimensional submanifold $S \subset M \backslash$ 
satisfies the condition $(\rm{Im} \rho_s)|_S \equiv 0$. According to Proposition 4 the differential of
this 1 - form equals to $2 \pi \omega$ which implies that $S$ must be lagrangian. The restrictions
of $a$ and $a_s$ on $S$ are both flat and their difference is real exact form $\rho_s|_S = d (\rm{ln} \vert s \vert|_S)$.
Since $s$ is covariantly constant w.r.t. $a_s$ it follows that $\frac{s}{\vert s \vert}|_S$
must be covariantly constant w.r.t. $a|_S$, therefore $S$ must be BS - submanifold. The speciality condition
obviously follows from the same argument, and it ends the proof.

Topology comes to the disccusion at this step: according to Proposition 4 our calibrating 1 -from $\rm{Im} \rho_s$
gives a cohomology class being restricted to any lagrangian submanifolds. Indeed, the restriction is a closed 1 - form,
and the corresponding cohomolgy class has been discussed above: 

{\bf Proposition 5.} {\it Let pair $(p(s), S) \subset \mathbb{P}(\Gamma(M, L)) \times {\cal B}_S$ belongs to a connected component
$K_i$ from the Proposition 1 above. Then the cohomology class $[\rm{Im}(\rho_s)] = m(K_i) \in H^1(S, \mathbb{Z})$.}   
  
Indeed, let $s|_S = \alpha(s, S) \sigma_S$ be as in Proposition 1. The proportionality coefficient 
$\alpha(s, S)$ is a complex valued non vanishing function on $S$. The logarithmic derivative of this function
is exactly our 1 - form  $\rho_s$ after restriction to $S$; the real part is exact so all cohomological data
is concentrated in $\rm{Im} \rho_S$, which gives the statement of Proposition  5.

SBS - submanifolds can be understood as zeros of certain vector field on a subset of the moduli space
${\cal B}_S$. Let us fix a smooth section $s \in \Gamma (M, L)$ and consider the connected components $K_i(s) \subset {\cal B}_S$, 
defined as $p_2(K_i \cap p_1^{-1} (p(s)))$ where $K_i$ were defined in Proposition 1 above,and then
take the components with the trivial classes  $m_i(K)$. Then for any point $[S] \in K_i(s)$ the restriction of the form
 $\rm{Im} \rho_s|_S$ is exact. But the exact forms are tangent vectors to  ${\cal B}$ at point $[S]$ (see [8], [9]), 
 therefore our section $s$ generates a smooth vector field $\tau_s$ on each appropriate component $K_i(s)$.
 The point is that zeros of this vector field $\tau_s$ present SBS - submanifolds. Moreover, this vector field
 $\tau_s$ is transversal: zeros of $\tau_s$ is isolated. And even more can be said for the field: recall that
 the moduli space ${\cal B}_S$ is covered by Darboux - Weinstein neighborhoods which play the role of charts
 in the natural atlas (see [8]). Then we have the following

{\bf Theorem 3.}  {\it Any Darboux - Weinstein neighborhood contains at most only one  $s$ - special Bohr - Sommerfeld
lagrangian submanifold w.r.t. a fixed smooth section $s$.}

 Let the ``center'' of a Darboux - Weinstein neighborhood --- a BS - submanifold $S_0 \subset M \backslash D_s$ ---
 is itself SBS w.r.t. to our fixed section  $s$. Then accroding to Theorem 2 the restriction  $\rm{Im} \rho_S|_S$ 
 identically vanishes.  Transport our 1- form  $\rm{Im} \rho_s$ to a small neighborhood of the zero section
 ${\cal O}_{\epsilon} (S_0) \subset T^* S_0$ and denote this form as $  \beta_s$. Then near zero section in 
 $T^* S_0$ we have two  1- forms  --- the canonical 1 -form $\alpha$ and new 1 - form $\frac{1}{2 \pi} beta_s$ 
 with the same differentials. The difference of these forms is closed and identically vanishes on the zero section.
 
 According to the Darboux - Weinstein theorem any close BS - submanifold of the same type is presented by the graph
 of an exact 1 -form  $df$ for certain real function $f \in C^{\infty}(S_0, \mathbb{R})$. Every function on
 a compact set admits ar least two critical points which we denote as $x_+, x_-$ for function $f$. These points 
 belongs to the intersection of the graph $\Gamma (df)$ and the zero section $S_0$  in $T^* S_0$.
 Consider two pathes  $\gamma_0, \gamma_f = \Gamma (df (\gamma_0))$, connecting the points along $S_0$ and $\Gamma(df)$ 
 respectively. Then the integrals along closed loop $\gamma_0 \cup \bar \gamma_f$ for the forms $\alpha$ and
 $\frac{1}{2 \pi}\beta_s$ must be the same since both the forms have the same differential. But if we suppose that
 the graph $\Gamma (df)$ is again special w.r.t. to the same section $s$, then the intergal for the second form
  $\frac{1}{2\pi} \beta_s$ must be trivial. Indeed, by the Theorem 2 our 1-form  $\beta_s$ must vanishes being restricted
  to  $\Gamma (df)$ as well as to the zero section. But the integral for $\alpha$ along the same loop
  equals to the difference  $f(x_+) - f(x_-)$, so it could happen iff the function $f$ is constant. But then
$\Gamma (df) = S_0$.

  The case when  $S_0$ itself is not special, but the graph of the differential of certain function
  $f_0 \in C^{\infty}(S_0, \mathbb{R})$ is special, is essentailly the same: if for some other $f$ 
  we again have special graph then  again a pair of points --- maximum and minimum of the difference $f - f_0$, --- 
  presents the intersection points for the graphs thus one again takes the pathes and comes to the contradiction.
  It ends the proof of Theorem 3.

The last theorem clarifies Proposition 2 above: the discretness stated there is of certain special type,
which depends on the ``values'' of Darboux - Weinstein neighborhoods for Bohr -Sommerfeld lagrangian submanifolds.
As we have seen in any such a neighborhhood one has two 1- forms: the canonical 1-form $\alpha$ and 
the form $\beta_s$ defined by our fixed smooth section $s$. Both of them have the differential equals to
the symplectic form. If the restriction $\beta_s|_{S_0}$ is exact then one can take a smooth function
$f_0 \in C^{\infty}(S_0, \mathbb{R})$ such that $df_0 = \beta_s|_{S_0}$. Then for each point $x \in S_0$
one can reconstruct a smooth function $\psi_x$ on the corresponding ball in $T_x^*S_0$ which is given
by the Darboux - Weinstein neigborhood near $x$, namely since the ball is contractible the restriction
of $\beta_s$ to it is exact and the integral $\int_x^p \beta_s$ depends on the end $p$ in the fiber
$T^*_x S_0$ only, so it gives as $\psi_x (p) = f_0 (x) + \int_x^p \beta_s$. Globalizing over $S_0$
we get function $\psi(x,p) = \psi_x(p)$, which by the construction satisfies
$$
\alpha - \beta_s = d \psi(x, p).
$$
Near the boundary of the Darboux - Weinstein neighborhood this function can be smoothly deformed to zero, and 
the denote the result again as $\psi$. Then if the fixed section $s$ is deformed by the following family
$s_t = s e^{-it\psi}, t \in [0, 1], s_0 = s,$ and we take the corresponding calibrating 1 - forms $\beta_{s_t} =
\beta_t$ then for each $t \in [0,1]$ one can define a non linear transformation of the space of exact 1 - forms
$B^1(S_0)$. For each $\beta_t$ take an exact 1 - form $\eta \in B^1(S_0)$ and consider 
$$
A_t(\eta) = \pi_* \beta_t|_{\Gamma (\eta)} \in B^1(S_0)
$$
where $\Gamma (\eta) \subset T^*S_0$ is the graph of $\eta$. Note however that this transformation is cerrectly
defined for small exact 1 - froms only since the graph $\Gamma(\eta)$ must lie in the Darboux - Weinstein
neigborhood.  This means that $A_t$ is defined for certain small ball in $B^1(S_0)$. Our claim is
that if $S_0$ is $s$ = $s_0$ special BS then $A_t$ is locally surjective for any $t$. Locality here means that
we are interested in a small ball around zero in $B^1(S_0)$. The claim is based on two facts: first,
for any $t \in [0,1]$ the map $A_t$ injective; second, for $t = 1$ the map is surjective being indentical
since for $t=1$ our form $\beta_1$ coincides with the canonical form $\alpha$. The proof of injectivity for
all $t$ follows the same arguments as for Theorem 3. Now let our fixed section $s$ is deformed to
$s_{\delta}$, and the corresponding 1 - form is $\beta_{\delta}$. Then we extend the family of
transformations $A_t, t \in [0,1]$ to $A_t, t \in [0, 1+ \delta]$ where $\delta$ is small enough.
Again all $A_t$ are injective, so one can expect that for sufficiently small $\delta$
they are localy surjective, and then for small $\delta$ zero 1 - form belongs to the image $A_{\delta}$
therefore the corresponding preimage gives Bohr - Sommerfeld deformation of $S_0$ which is $s_{\delta}$-
special. This ends the sketch of the proof of

{\bf Theorem  4.}{\it The image of the first projection $\rm{Im} p_1 \subset \mathbb{P} (\Gamma (M, L))$ is open subset.}

Combinig Theorem 1 and Theorem 4 we get

{\bf Corollary.} {\it The space ${\cal U}_{SBS}$ admits a Kahler structure.}

   Indeed, we can just lift the standard Kahler structure from $\mathbb{P} (\Gamma (M, L))$ using $p_1$.

Thus the space ${\cal U}_{SBS}$ can be regarded as certain ``complexification'' of
the moduli space ${\cal B}_S$, however it is not a complexification in usual sense: 
the ``dimension'' of  ${\cal B}_S$ is much less that the ``half dimension'' of  ${\cal U}_{SBS}$.  

\section{Algebraic case}

Now suppose that our symplectic manifold $(M, \omega)$ admits a compatible complex structure $I$ which is integrable.
This means that $M$ is endowed with a Kahler metric of the Hodge type (since $\omega$ defines an integer cohomology class)
therefore $(M, \omega, I)$ can be regarded as an algebraic variety, see [10]. On the other hand the prequantization data
$(L, a)$ in the case correspond to a holomorphic line bundle since the curvature $F_a$ has type (1,1) w.r.t. $I$
therefore it induces a holomorphic structure on our hermitian line bundle $L$. This means that fixing $I$ we 
cut a finite dimensional subspace $H^0(M_I, L) \subset \Gamma (M, L)$ formed by holomorphic sections of $L$.
This subspace is finite dimensional, and it is natural to construct the corresponding finite dimensional vesrion
of the space ${cal U}_{SBS}$ defined above. Take in the direct product $\mathbb{P} H^0(M_I,L) \times {\cal B}_S$
the subset defined by the specialty condition as it was above and get {\it the moduli space} of SBS lagrangian submanifolds
over $(M, \omega, I)$ defined as the preimage
$$
{\cal M}_{SBS} = p_1^{-1} (\mathbb{P} H^0 (M_I, L)) \subset {\cal U}_{SBS}
$$
of the first projection of the projectivization of the holomorphic section  subspace.

Since the propreties of projection $p_1$ were studied above we know that
$$
p_1: {\cal M}_{SBS} \to \mathbb{P}H^0(M_I, L)
$$
has discrete fibers, so the moduli space ${\cal M}_{SBS}$ is a finite dimensional set fibered over an open subset in
the projective space $\mathbb{P}H^0(M_I, L)$. In the rest of the present text we show that for a generic
holomorphic section the number of preimages at the corresponding fiber is finite and propose a constructive way
how to find SBS lagrangian submanifolds. On the other hand in the known examples the open part of $\mathbb{P}H^0(M_I, L)$
has very natural form: we just remove an algebraic subvariety from the projective space. This hints that the moduli space
${\cal M}_{SBS}$ admits certain natural compactification.

The key observation for the algebraic case is the following:

{\bf Proposition 6.} {\it For a holomorphic section $s \in H^0(M_I, L)$
the corresponding calibrating 1- form $\rm{Im}(\rho_s)$  equals to $- I(d(\rm{ln} \vert s \vert)$.}

First, the form $\rho_s$ has type (1,0) w.r.t. the complex structure $I$. Indeed, since
$$
\rho_s = \frac{<\nabla_a s, s>}{<s,s>}
$$
but $s$ is holomorphic w.r.t. $\partial_a$ therefore $\nabla_a s \in L \otimes \Omega^{1,0}$ and the resting
operations do not change the type.

Second, for the real and imaginary parts of a (1, 0) - form we have the corresponding relation, and since we know the real part
of $\rho_s$ which equals to $d (\rm{ln} \vert s \vert)$ we get the statement of Proposition 6.

Therefore for a holomorphic section $s \in H^0(M_I, L)$ we get the corresponding real smooth function
$$
\phi_s = - \frac{1}{2\pi}\rm{ln} \vert s \vert
$$
which is correctly defined on the complement $M \backslash D_s$ and which is {\it plurisubharmonic} 
or {strongly convex} w.r.t. $I$ on the complement
since $d (I (d \phi_s)) = d^c d \phi_s = \omega$ implied from  Proposition 4. The convexity property is very usefull
for our investigations since we have the following remark

{\bf Proposition 7.}{ \it A lagrangian submanifold $S \subset M$ is SBS w.r.t. a holomorphic section $s \in H^0(M_I, L)$
iff $\rm{grad} \phi_s$ is parallel to $TS$ at each point of $S$.}

Indeed, we know from Theorem 2 above that SBS condition is equivalent to $\rm{Im} \rho_s|_S \equiv 0$, but
if $S$ is lagrangian and $\rm{Im} \rho_s = -I (d \rm{ln} \vert s \vert)$ then the SBS condition is equivalent to
the fact that $\omega (v, \rm{grad} \phi_s)$ identically vanishes if $v$ is tangent to $S$ at each point.
Thus $\rm{grad} \phi_s$ must be parallel to $TS$.

Now suppose that for a holomorphic section $s \in H^0(M_I, L)$ the corresponding function $\phi_s$ is a Morse function
on the complement $M \backslash D_s$. It means that the real function $\vert s \vert^2$ is {\it relatively} Morse
so it admits degenerated absolute minimum at  $D_s$ but all other critical points are non degenerated. For this case
consider the critical points $x_1, ..., x_k \in M \backslash D_s$ of the function $\phi_s$ .
This function being restricted to any  compact lagrangian submanifold $S$ must have critical points (not less that
the Morse inequality dictates) and if $S$ is SBS w.r.t. $s$ it implies that every critical point of the restriction
$\phi_s|_S$ must be critical for global $\phi_s$ at the same point of $S$ considering as a point of $M \backslash D_s$.
Indeed, at these points both $d \phi_s$ and $I (d \phi_s)$ vanish being restricted to $S$ (note that we suppose that
$S$ is SBS so the last form vanishes identically on $S$). Consequently any SBS lagrangian submanifold must
contain several critical points from the set $x_1, ..., x_k$. For example if $S$ is a lagrangian torus then
the number of critical points it must contain is not less than $n$. And if for a generic holomorphic section
$k< n$ (as it is for the case of $\mathbb{C}\mathbb{P}^2$ and $L = {\cal O}(2)$) then SBS lagrangian tori do not exist.

Moreover, Proposition 7 shows that the gradient flow generated by $\phi_s$ must preserve SBS lagrangian submanifold
$S$. It means that SBS submanifold $S$ must contain not just critical points but trajectories of the gradient flow.
In simple cases it hints how to completely solve the problem.

{\bf Example.} Consider the projective line $\mathbb{C} \mathbb{P}^1$ with the standard Kahler structure and rescale
the Kahler form by 2. The corresponding line bundle is ${\cal O}(2)$ so any holomorophic section
$s$ is completely determined by its zeros up to $\mathbb{C}^*$, and since our SBS condition depends on
the class up to $\mathbb{C}^*$ it means that every pair of points $(x_1, x_2)$ should define SBS lagrangian loops.
Suppose first that the section is multiple so the points coincide $x_1 = x_2$. Then the function
$\phi_s$ has only one critical point $x_{min}$ on $\mathbb{C} \mathbb{P}^1$ and for this section
SBS lagrangian loop doesn't exist. Indeed, if it exists then $\phi_s|_S$ must have at least two different
critical points which must be global critical points of $\phi_s$ on the punctured sphere but there
one has only one critical point $x_{min}$. In contrast for a section $s$ with two different zeros
$x_1, x_2$ a SBS loop must exist due to pure topological reasons. Indeed, the vector field $\rm{grad}
\phi_s$ in this case has at least three singular points $x_1, x_2, x_{min}$ with positive indecies.
But the Euler characteristics of sphere is 2 which implies that $\phi_s$ admits a saddle point $x_s$.
Therefore the gradient flow must have one separatrix $\gamma$ passing $x_{min}$ and $x_s$
which separates trajectories which go to $x_1$ and $x_2$. This $\gamma$ is preserved by the gradient flow
which implies that $\gamma$ is SBS with respect to $s$. There are exceptional cases when
the zeros of $s$ are antipodal, and in this case $\phi_s$ is not Morse on the complement but
admits a critical equatorial loop which is again SBS. 

This presents another principle detecting
SBS submanifolds: if for a holomoprhic $s$ the function $\phi_s$ is not Morse but has sufficiently
big critical subset (of maximal possible dimension $n$) it must be SBS  lagrangian w.r.t. this section.
For example for $\mathbb{C}\mathbb{P}^2$ and $L = {\cal O}(2)$ the holomorphic section $s$ corresponding
to conic $\sum z_i^2 = 0$ defines $\phi_s$ which is not Morse on the complement but whose critical set
is exactly the real part $\mathbb{R} \mathbb{P}^2 \subset \mathbb{C}\mathbb{P}^2$ of the projective plane
which gives us an example of SBS lagrangian $\mathbb{R}\mathbb{P}^2$.

All these show that the set $\rm{Crit} \phi_s \subset M \backslash D_s$ of critical for $\phi_s$ points
together with the lines of the gradient flow of $\phi_s$ connecting the critical points of ``finite type''
(so we exclude infinite maximums of $\phi_s$ at zeros $D_s$) 
form  in a sense the ``base set'' $B$ for SBS lagrangian submanifolds. More precisely,  consider 
on the complement $M \backslash D_s$ the set of critical points of $\phi_s$. 

{\bf Definition 3.} {\it For a holomorphic section $s \in H^0(M_I, L)$ we define the base set
$B_s \subset M \backslash D_s$ as a subset of the complement $M \backslash D_s$
which contains all critical points of $\phi_s$ together with all finite trajectories
of the gradient flow of $\phi_s$ connecting pairs of critical points. }

So if a pair $x_k, x_l$ of 
the critical points is connected by a finite trajectory $\gamma$ then the points of $\gamma$
together with $x_k$ and $x_l$ must belong to $B$.

Proposition 7 implies the following corollary: every SBS lagrangian submanifold $S$ must lie
in $B$. 

According to an old result of Milnor, [12], if $\phi_s$ is strongly convex on the complex manifold
$M_I \backslash D_s$ then every Morse critical point $x_i$ of $\phi_s$ has the Morse index less or equal
to $n$ where $n$ is the complex dimension of $M_I$. Moreover, the negative subspace $T^-_x M_I \subset T_x M_I$
which correspond to incoming trajectories for the gradient flow must be isotropical w.r.t. our
Kahler form $\omega$. Suppose that $S \subset M \backslash D_s$ is special  w.r.t. $s$
and suppose that $S$ is compact. The restriction of $\phi_s$ to $S$ admits at least one maximum $x_m \subset S$.
Due to the arguments around Proposition 7 we know that $x_m$ is a critical point of $\phi_s$ on $M \backslash
D_s$. But it implies that $x_m$ must have the Morse index equals exactly $n$. The corresponding $n$ -dimensional family
of incoming trajectories of the gradient flow must be contained by $S$. This implies the fact that the number
of SBS lagrangian submanifolds is less or equal to the number of critical points of $\phi_s$ os the Morse index
$n$. Thus we have sketched the proof of the following 

{\bf Theorem 5.}{\it For a generic holomorphic section $s \in H^0(M_I, L)$ the number of SBS lagrangian
submanifolds is finite.}

At the same time the restriction on the number of SBS submanifolds doesn't automatically imply the existence
of a single one, but the Morse theory for strongly convex functions $\phi_s$ says that the existence theorem
can be formulated after certain extension of the definition of SBS - submanifolds to SBS -cycles.
This extension can be illustrated by the following 

{\bf Example.} Consider $M = \mathbb{C} \mathbb{P}^1$ and $L = {\cal O}(d)$. In this case holomorphic
sections of $L$ are presented by homogenious polynomials of degree $d$ in variables $[z_0: z_1]$.
In [13] one finds the separatrix trajectories for the gradient flow generated by generic golomorphic sections
of ${\cal O}(d)$ using the natural reformulation of the problem in terms of the polynomials. Taking
the section $s \in H^0(\mathbb{C} \mathbb{P}^1, {\cal O}(d)$ which corresponds to $P_d = z_0^d + z_1^d, d> 2,$ 
we get for the function $\phi_s$ three types of singular points on $\mathbb{C} \mathbb{P}^1$:
first, two minimal points at $[1:0]$ and $[0:1]$; second, $d$ saddle points at $[1: e^{2 \pi k i/d}]$
and, third, $d$ infinite maximal points coming with zeros of $P_d$. Thus as the separtrix trajectories
  we have $d$ open segments $\gamma_k, i =1, ..., d,$ 
 joining the poles and passing through the corresponding saddle point $[1: e^{2\pi k i/d}]$.
 In the even case we have the coincidence of tangents to $\gamma_k $ and $\gamma_{k + d/2}$
 at the poles $[1:0], [0:1]$ therefore we can combine closed smooth separatrix trajectories
 but it doesn't happen if $d$ is odd. Moreover, it is 
  not hard to see that for a generic polynomial of degree $d$ we should get essentially
 the same picture, but losing the symmetry of $P_d$ even in the case of even $d$ we should lost 
 smooth combined trajectories: the coincidences of the tangents at the poles should be lost.
 This hints the way how the theory can be regularized: we must allow submanifolds with certain
 types of singularities and extend the considerations from SBS submanifolds to
 SBS cycles. In the case $\mathbb{C} \mathbb{P}^1, {\cal O}(d)$
 we allow loops with finite numbers of ``corners'' and it leads to the following description:
 there are $d(d-1)/2$ SBS cycles for a generic holomorphic section of ${\cal O}(d)$.
 In [13] one studies in details other fibers of the projection $p_1: {\cal M}_{SBS} \to \mathbb{P}
 H^0(\mathbb{C} \mathbb{P}^1, {\cal O}(d))$ and describes the ramification structure of it.
 
 Thus we shall consider not only lagrangian embedding but as well lagrangian immersions $S \subset M$,
 and, follow [8], we call these $S$ lagrangian cycles. As it was pointed out in [8] the Bohr -
 Sommerfeld condition can be imposed on lagrangian immersions as well as on smooth embedding:
 it is clear from the remarks after Definition 1 above since for immersed $S \subset M$ one can consider
 loops on $S$, discs with boundaries on the loops and the symplectic area of the discs which is an integral
 so it is correctly defined even in the case of non smooth loops. 
 
 The speciality condition can be extended to the immersions as well: we say that a BS lagrangian cycle
 is special w.r.t.  a holomorphic section $s \in H^0(M_I, L)$ iff the calibration form
 $\rm{Im} \rho_s = - I (d \rm{ln} \vert s \vert)$ identically vanishes on $S$. 
 For a singular point of $S$ it means that any tangent vector from the tangent cone annihilates the calibration form. 
 As we have seen it could happen if the singular point of $S$ is a critical point of the Kahler potential $\phi_s$.
 
 The main conjectures which can be formulated at the present time are the following
 
 {\bf Conjecture.} {\it The number of special Bohr - Sommerfeld lagrangian cycles is invariant for
 generic holomorphic sections from $H^0(M_I, L)$.}
 
 If this is true then one gets certain system of lagrangian invariants for algebraic varieties. 
 
 On the other hand as it was  mentioned above we are interested in the moduli spaces ${\cal M}_{SBS}$
 which now is defined as the moduli space of SBS lagrangian {\it cycles}. These moduli spaces admitsthe same as above projections
 to the projectivized spaces of holomorphic sections thus we can expect that these moduli spaces are finite dimensional
 Kahler manifolds which can be naturally compactified. This is our second conjecture. 
 
 \section{Example: complex 2 -dimensional quadric}

  In this section we discuss the first ``non toy'' example and present some technical arguments which should be
  exploited in general case of algenraic varieties.
  
  Our $M$ is complex 2 - dimensional quadric $Q$ relized either as the direct product $\mathbb{C} \mathbb{P}^1 
  \times \mathbb{C} \mathbb{P}^1$ or as a subvariety in $\mathbb{C} \mathbb{P}^3$ defined by a quadratic 
  polynomial. Our line bundle $L$ is taken to be ${\cal O}(1,1)$ so the tensor product of
  two copies of ${\cal O}(1)$ on both $\mathbb{C} \mathbb{P}^1$'s lifted by the projections to the
  direct summands. The symplectic form is taken to be the direct sum of the lifted from the direct summands
  Fubiny - Study forms. All details on Kahler geomtry can be found in [10].
  
  Recall first that $Q$ admits lagrangian 2 - spheres: to see this let us fix coordinates $[x_0: x_1]$
  and $[y_0: y_1]$ on the projective lines and consider the subset
  $$
  S_0 = \{ y_i = \bar x_i\} \subset \mathbb{C} \mathbb{P}^1 \times \mathbb{C} \mathbb{P}^1
  $$
  which is obviously 2 - dimensional sphere, and short calculation shows that it is lagrangian.
  Note that every lagrangian sphere automatically satisfies the Bohr - Sommerfeld condition.
  
 Then holomorphic sections of ${\cal O}(1,1)$  are presented by linear combinations $\sum \alpha_{ij} x_i y_j$
  and it is not hard to see that $S_0$ is special Bohr - Sommerfeld w.r.t. the section given by
  $s_0 = \{\alpha_{00} = \alpha_{11} = 1, \alpha_{01} = \alpha_{1,0} =0\}$. Are there other $s_0$- SBS lagrangian spheres?
  The answer is negative: we take the function $\phi_{s_0}$ and find that the base set $B_s \subset Q \backslash D_{s_0}$
  from Definition 3 above in this case   coincides with $S_0$ so no other SBS spheres. 
  
  If we realize $Q$ as a quadratic surface in $\mathbb{C} \mathbb{P}^3$ then the holomorphic sections
  of ${\cal O}(1,1)$ up to scale are presented by projective planes $H_s \subset \mathbb{C} \mathbb{P}^3$ 
  so the zero sets $D_s$ are presented
  by the intersections  $H_s \cap Q$, see [10]. We have distingiushed sections given by
  tangent planes so the intersections and consequently $D_s \subset Q$ for these cases are not smooth.
  For generic $s$ the subvariety $D_s$ is a smooth conic topologically equivalent to $S^2$, and for
  the tangent planes nd the corresponding sections $D_s$ are given by pair of intersecting projective lines
  so topologically these are equivalent  to pair of $S^2$ tranversally intersecting at a point.
  Such a section is given by a reducible expression $\sum \alpha_{ij} x_i y_j =
  (a_0x_0 +a_1 x_1)(b_0 y_0 + b_1 y_1)$, and these two lines are $[-a_1: a_0] \times \mathbb{C} \mathbb{P}^1
  \cup \mathbb{C} \mathbb{P}^1 \times [-b_1: b_0]$. Since the hermitian norm of a holomorphic section
  of ${\cal }(1)$ on $\mathbb{C} \mathbb{P}^1$ has exactly two critical points it is not hard to
  see that the tensor product of two sections lifted to $Q$ from the direct summands must have
  exactly one critical point on $Q \backslash D_s$ therefore the base set $B_s$ for a reducible
  section consists of exactly one point. This implies that for reducible holomorphic sections
  there are no SBS lagrngian cycles.
  
  Now we claim that for a generic holomorphic section $s \in H^0(Q, {\cal O}(1,1)$
  there exists unique SBS lagrangian sphere. Indeed, if $s$ is a generic holomorphic section
  then the function $\phi_s$ on $Q \backslash D_s$ has exactly two critical points:
  minimal point $p_{min}$ and ``saddle'' point $p_s$ with Morse index 2. Then the base set $B_s$ is again
  a 2 -sphere: we take two dimensional family of incoming to $p_s$ trajectories whose tangents at $p_s$
  span  the negative definite subspace for the Hessian. Due to the Milnor remark ([12]), this subspace must be lagrangian,
  but the Liouville vector field $grad \phi_s$ preserves the lagrangian condition, so reversing the gradient
  flow of $\phi_s$ we get that the tangent space to $B_s$ at each point is lagrangian. Therefore we get a lagrangian sphere
  which is by the construction SBS w.r.t. the choosen section.
  
  Thus the main strategy in construction of SBS cycles is the Morse theory of $\phi_s = - \rm{ln} \vert s \vert$
  for holomorphic sections! It is related to the Morse theory of $\vert s \vert^2$ which is never Morse
  (we have $D_s$ as the subset of non isolated minima) but the case when $\vert s \vert^2$ has
  isolated non degenrated critical points on the complement $M \backslash D_s$ can be studied.
  Here we present the arguments for $Q$ when $D_s$ is smooth. Suppose that $s$ is {\it relatively Morse}
  so that $\vert s \vert^2$ is Morse on the complement $Q \backslash D_s$. The gradient of $\vert s \vert^2$
  is outgoing near $D_s$: for any small tubular neigborhood of $D_s$ at the boundary points the gradient vectors see
  outside of the neigborhood. It is well known that the subset of Morse functions is dense in the function space
  therefore $\vert s \vert$ can be deformed to a Morse function $\psi$ which has the same critical points
  on $Q \backslash D_s$. But let us choose $\psi$ in such a way that $\psi|_{D_s}$ is strictly Morse
  so it has only two critical points --- maximal and minimal. Due to the fact that the gradient is outgoing
  near $D_s$ it follows that $\psi$ has exactly two additional critical points of indecies 0 (the minimum on
  $D_s$) and 2 (the maximum on $D_s$) to the set of critical points of $\vert s \vert$. If $\psi$ is strictly Morse
  we are done: the topology of $Q$ dictates that it must be 4 critical points of indecies 0, 2, 2, 4, and the last
  pair is the desired $p_{min}$ and $p_s$ for our base set $B_s$. Another pair of related by the gradient flow
  critical points with the index difference 2 cann't exist - otherwise we would get another cell in $H^2(Q, \mathbb{Z}$
  (it is interesting that in this picture two cells from $H^2(Q, \mathbb{Z})$ are realized by holomorphic ($D_s$)
  and lagrangian 2 - spheres.
  
  What happens if the section $s$ is reducible? Then $D_s$ takes not 2 but 3 critical points of $\psi$
 and for $\phi_s$ on $Q \backslash D_s$ it remains unique critical point.
 
 All these facts can be checked directly via computations for critical points of $\vert s \vert^2$.
 in coordinates $[x_0: x_1], [y_0:y_1]$: we take the exrpession 
 $$
 F_{\lambda, \mu} = \vert \sum \alpha_{ij}x_iy_j \vert^2 - \lambda(\vert x_0 \vert^2 + \vert x_1 \vert^2 -1)
 - \mu (\vert y_0 \vert^2 + \vert y_1 \vert^2 -1)
 $$
 and solve the system
 $$
 \frac{\partial F_{\lambda, \mu}}{\partial z} = 0
 $$
 where $z = x_i, \bar x_i, y_j, \bar y_j$ (as usual for the computations of conditional extremums). 
 The system can be solved which gives us the set of critical points for the section $s$.
 
 Summing up we get the following answer: for the case $M = Q$
 with the standard Kahler structure and $L = {\cal O}(1,1)$ with the standard hermitian structure,
 for the topological data $S \cong S^2$, $[S] = (1, -1) \in H^2(Q, \mathbb{Z})$, the moduli space
 ${\cal M}_{SBS}$ is naturally isomorphic to $\mathbb{C} \mathbb{P}^3 \backslash Q'$
 where $Q'$ is a quadric (projetively dual to $Q$). Another fact: this moduli space has
 a natural compactification isomorphic to $\mathbb{C} \mathbb{P}^3$.

  \section{Final remark}

At the end we would like to mention the following natural construction.

In [14] one constructs a natural $\mathbb{C}$ - bundle ${\cal L}$ with a fixed hermitian
structure over the moduli space ${\cal B}_S$. On the other hand the projective space  
 $\mathbb{P} H^0(M_I, L)$ is naturally endowed with the line bundle ${\cal 0}(1)$ together with
 a hermitian connection  $a$, whose curvature form is proportional to the Kahler form.
 In the presence of a complex structure $I$ on $M$ the bundle ${\cal L}$ 
 is naturally endowed with a connection $A_I$ which depends on the choice of $I$. Therefore 
 on the moduli space ${\cal M}_{SBS} \subset \mathbb{P} H^0(M_I, L) \times {\cal B}_S$ 
 of SBS lagrangian cycles one naturally gets a bundle  ${\cal E}$, which is the tensor product 
$p_1^* {\cal O}(1) \otimes p_2^* {\cal L}$ restricted to the moduli space ${\cal M}_{SBS}$. This bundle is endowed
with a hermitian connection $\mathbb{A}$, given by the connections $a$ and $A_I$ on the tensor product components. 
Since the cuvrature of such a connection equals to the sum of the curvatures of $a$ and $A_I$ if the last
one is trivial then it were possible to consider $\mathbb{A}$ as a connection whose curvature
is again proprotional to the Kahler form (recall that we lift the Kahler structure
from the projective space). Then we would come to the situation when all is ready
for further ``quantization''...

Recall the construction from [14]: in the direct product $M \times {\cal B}_S$ one has the ``incidence cycle'' 
$$
{\cal N} = \{ (x, S) \quad | \quad x \in S \} \subset M \times {\cal B}_S
$$
with natural projections $q_i$ to the direct summands. Then the lift $q_1^* (L, a)$ admits one - dimensional
0- cohomology space along the fibers of $q_2$ (this follows from the Bohr - Sommerfeld condition),
therefore
$$
{\cal L} = R^0 (q_2)_* q_1^* (L, a) \to {\cal B}_S
$$
is a line bundle.   Over point $S \subset {\cal B}_S$ the fiber is spanned by the section $\sigma_S$, 
consequently there is the natural $U(1)$ -action, generated  by the $U(1)$ -action
on the prequantization bundle $L$: element $e^{i c} \in U(1)$ just twists all  $\sigma_S \mapsto e^{ic} \sigma_S$.

The connection $A_I$ on this bundle ${\cal L} \to {\cal B}_S$ can be constructed as follows.
For any point $S \subset {\cal B}_S$ consider the Darboux - Weinstein neigborhood
${\cal O}_{DW}(S) \subset M$, so all close to  $S$ points of the moduli space ${\cal B}_S$ are given by
graphs $S_f = \Gamma (df)$ of the differential for functions $f \in C^{\infty}(S, \mathbb{R}$,
and the corresponding covariantly constant sections of the restrictions  $(L,a)|_{S_f}$ are given by twisting
of the form $e^{i(f+c)} \sigma_S$. If there is a natural way how to choose for any function $f$ a representative from
the class $ f+ c, c \in \mathbb{R}$ universally then it gives a local section of the bundle ${\cal L}$ over the
neigborhood $S \subset {\cal B}_S$. In general this universal choice doesn't exist, however in our case
when a compatible complex structure on $M$ is fixed therefore we have a fixed riemannian metric $g$, 
then we can take the restriction $g|_S$ and define the corresponding volume form $d \mu(g|_S)$. Then
the natural condition $\int_S f d \mu(g|_S) = 0$ specifies an $f$ from the corresponding class $C^{\infty}(S, \mathbb{R})/const$
(and this choice doesn\t depend on the orientation of $S$). Since for $S$ and $S_f$ the restrictions of $g$ {\it a priori}
give different volume forms it doesn't directly lead to a local section of ${\cal L}$
but for each point it gives a horizontal subspace in the corresponding point of the tangent bundle $T (\rm{tot} {\cal L})$, 
and the universal $U(1)$ - action lifts this horizonal subspace to a horizontal distribution
which is by the construction $U(1)$ - invariant. This is our hermitian connection $A_I$.

Now the problem is to calculate the curvature for this connection $A_I$. 
If it is flat then we would get on the moduli space ${\cal M}_{SBS}$ a hermitian line bundle
together with a hermitian connection whose curvature form is proportional to the Kahler form. 

The work on the problems listed above is in progress, so one hopes that new results will be found 
soon in special Bohr - Sommerfeld geometry.

{\bf References:}

[1] {\bf Yu. I. Manin}, {\it ``Foreword for the 3d volume''}, A.N. Tyurin's selected papers, (in Russian),
Moscow - Izhevsk,  2004; 

[2] {\bf M. Kontsevich}, {\it ``Homological algebra of mirror symmetry''}, Proceedings of ICM (Zurich, 1994),
Birkhouser, Basel 1995, pp. 120 - 139;

[3] {\bf A.N. Tyurin}, {\it ``Geometric quantization and mirror symmetry''}, arXiv: math/9902027v1;

[4] {\bf N. Hitchin}, {\it ``Lectures on special lagrangian submanifolds''},  Winter school on Mirror 
Symmetry (Cambridge MA, 1999), AMS/IP Stud. Adv. Math,
23, AMS 2001, Providence,  pp. 151 - 182;

[5] {\bf A. Strominger, S.-T. Yau, E. Zaslow}, {\it ``Mirror symmetry is T - duality''}, Winter school on Mirror 
Symmetry (Cambridge MA, 1999), AMS/IP Stud. Adv. Math,
23, AMS 2001, Providence,  pp. 333 - 347;

[6] {\bf D. Auroux}, {\it ``Mirror symmetry and T - duality in the complement of an anticanonical divisor''},
J. Gokova Geom. Topol. (2007), pp. 51 - 91;

[7] {\bf N.A. Tyurin}, {\it ``Special lagrangian fibrations on the flag variety''}, Theor. and Math. Phys.
167: 2 (2011), pp. 193 - 205;

[8] {\bf A.L. Gorodentsev, A.N. Tyurin}, {\it ``Abelian lagrangian algebraic geometry''}, Izvestiya math. 65: 3
(2001), pp. 15 - 50;

[9] {\bf N.A. Tyurin}, {\it ``Geometric quantization and algebraic Lagrangian geometry''}, LMS, Lecture
Note series 338, Cambridge 2007, pp. 279 - 318;

[10] {\bf P. Griffits, J. Harris}, {\it ``Principles of Algebraic geometry''}, NY, Wiley, 1978;

[11] {\bf R. Harvey, H. Lawson}, {\it ``Calibrated geometries}, Acta Math., 148 (1982), pp. 47 - 157;

[12] {\bf Y. Eliashberg}, {\it ''Topological characterization of Stein manifolds of $dim > 2$``},
Internat. J. Math., 1, no. 1, pp.
29 -46 (1990);

[13] {\bf A.N. Tyurina}, {\it ''Special Bohr - Sommerfeld lagrangian cycles on $\mathbb{C} \mathbb{P}^1`$``}, Batchelor diploma,
NRU HSE (Moscow, 2016), in preparation;

[14] {\bf N.A. Tyurin}, {\it  ``Algebraic lagrangian geometry: three geometric observations''},
Izvestiya math. 69 : 1 (2005), pp. 179 - 194. 

$$$$
Uspeniev den', Dubna
\end{document}